\documentclass[12pt,reqno]{amsart}
\usepackage{amsfonts}
\usepackage{indentfirst}
\usepackage{latexsym,amsmath,amsthm,amssymb,amsfonts}
\usepackage{setspace}
\usepackage{color}
\usepackage{tikz}
\numberwithin{equation}{section}
\allowdisplaybreaks

\def\CS{\mathcal S}

\def\vol{\mbox{vol}}
\def\s{\,\,\,\,}
\def\R{\mathbb{R}}

\def\C{\mathbb{C}}

\def\endproof{$\hfill\Box$\\}
\def\s{\,\,\,\,}
\def\R{\mathbb{R}}
\def\C{\mathbb{C}}

\numberwithin{equation}{section}
\newtheorem{theorem}{Theorem}[section]
\newtheorem{lem}[theorem]{Lemma}
\newtheorem{thm}[theorem]{Theorem}

\newtheorem{cor}[theorem]{Corollary}
\newtheorem{defi}[theorem]{Definition}
\newtheorem{rem}[theorem]{Remark}

\newcounter{Cnumber}

\title[John-Nirenberg Radius and Collapsing in
Conformal Geometry]
{\bf John-Nirenberg Radius and Collapsing in Conformal Geometry}
\author[Y. Li, G. Wei, Z. Zhou]
{Yuxiang Li,\, Guodong Wei,\, Zhipeng Zhou}
\address{\newline
Yuxiang Li:
 Department of Mathematical Sciences, Tsinghua University, Beijing 100084, P.R. China.
{\tt Email:yxli@math.tsinghua.edu.cn}
\newline
\newline
Guodong Wei:
Academy of Mathematics and Systems Science, CAS,
Beijing 100190, P.R. China.
{\tt Email:weiguodong@amss.ac.cn}
\newline
\newline
 Zhipeng Zhou:
Academy of Mathematics and Systems Science, CAS,
Beijing 100190, P.R. China.
{\tt Email:zhouzhipeng113@mails.ucas.ac.cn}}
\date{}
\begin{document}
\maketitle

\begin{abstract}
Given a positive function $u\in W^{1,n}$, we define its  John-Nirenberg radius at point $x$ to be the
supreme of the radius such that $\int_{B_t(x)}|\nabla\log u|^n<\epsilon_0^n$ when
$n>2$, and $\int_{B_t(x)}|\nabla u|^2<\epsilon_0^2$  when $n=2$. We will show that for a collapsing sequence of metrics in a fixed conformal class under some curvature conditions, the radius is bounded below by a positive constant. As applications, we will study the convergence of a conformal metric sequence on a $4$-manifold with bounded $\|K\|_{W^{1,2}}$, and prove a generalized H\'elein's Convergence Theorem.
\end{abstract}

\section{Introduction}
We say that a Riemannian manifold sequence collapses, if it converges to
a low dimensional space in the Gromov-Hausdorff distance. When $(M_k,g_k)$
collapses, a reasonable attempt is to  blow up the sequence, i.e., to find $c_k\rightarrow+\infty$, such that $(M_k,c_kg_k)$ converges to a manifold of the same dimension. This usually needs some  monotone properties,  such as volume comparison. Then some  sectional or  Ricci curvature conditions
are usually assumed for a collapsing sequence.

Recently, in \cite{Li-Zhou} the first author and the third author of this paper considered collapsing sequences in a fixed conformal class with bounded $L^p$-norm of scalar curvature, where $p>\frac{n}{2}$. Let $B_1$ be the unit ball of $\R^n$ centered at the origin and $g$ be a smooth metric over $\bar{B}_1$, where $n>3$. Consider a sequence of metric $g_k=u_k^\frac{4}{n-2}g$ which
satisfies
$$
\int_{B_1}|R(g_k)|^pdV_{g_k}<\Lambda,
$$
where $R(g_k)$ is the scalar curvature of $g_k$. Our conclusion is the following: ``when $\vol(g_k)\rightarrow 0$, there exists a sequence $\{c_k\}$ which tends to $+\infty$, such that $c_ku_k$ converges to a positive function in $W^{2,p}$ weakly". The proof of the conclusion is rather analytic and the John-Nirenberg inequality plays an essential role in the procedure.

Recall that the John-Nirenberg inequality says that, given
$u\in W^{1,q}(B_1)$, where $q\in[1,n]$ and $B_1$ is the unit ball of $\R^n$, if
$$
\int_{B_r(x)}|\nabla u|^q<r^{n-q},\s\s \forall \ B_r(x)\subset B_1,
$$
then there exists $\alpha$ and $\beta$, such that
$$
\int_{B_1}e^{\alpha u}\int_{B_1}e^{-\alpha u}<\beta.
$$

Inspired by the John-Nirenberg inequality we define the John-Nirenberg radius of $u_k$ at $x$ as follows:
\begin{equation}\label{radius1}
\rho_k(x)=\sup\left\{r:t^{2-n}\int_{B_t(x)}|\nabla\log u_k|^2<\epsilon_0^2,\s
\forall \ t<r\right\}.
\end{equation}
The key ingredient of the arguments in \cite{Li-Zhou} is that, when $\vol(g_k)$ converges to $0$, there must exist an $a>0$ which
is independent of $u_k$, such that $\inf_{B_\frac{1}{2}}\rho_k(x)>a$. This means that $$t^{2-n}\int_{B_t(x)}
|\nabla \log u_k|^2<\epsilon_0^2$$ for any $t<a$ and $x\in B_\frac{1}{2}$, hence
the John-Nirenberg inequality holds for $\frac{\log c_ku_k}{\epsilon_0}$ on $B_a(x)$, where $\int_{B_{\frac{1}{2}}}\log c_ku_k=0$. Then it follows the estimates of $L^\frac{\alpha}{\epsilon_0}$-norms of $\frac{1}{c_ku_k}$ and $c_ku_k$.

The arguments and calculations of the first half of \cite{Li-Zhou} were so complicated that it is not easy for one to pay attention to the John-Nirenberg radius, which was introduced and discussed in the last section of \cite{Li-Zhou}. While we think this new technique is very interesting and believe that it might be applied to some other nonlinear equations, we write this paper to highlight on the John-Nirenberg radius and give a simple explanation of how the John-Nirenberg inequality works.

It is easy to check that if  $$t^{2-n}\int_{B_t(x)}|\nabla u|^2<\epsilon_0^2$$ in \eqref{radius1} is replaced by $$t^{q-n}\int_{B_t(x)}|\nabla u|^q<\epsilon_0^q,$$ the arguments in
\cite{Li-Zhou} still work. We discover that it is much more convenient to use $q=n$ to define the John-Nirenberg radius. For this situation, the John-Nirenberg inequality can be
deduced from Moser-Trudinger inequality, which also gives the optimal constant in the John-Nirenberg inequality in the case of $q=n$. So we start our discussion
from Moser-Trudinger inequality in Section 2, and define the John-Nirenberg radius to be the supreme of the  radius such that
$$\int_{B_t(x)}|\nabla \log u|^ndx<\epsilon^n_0.$$
Then, we prove Theorem \ref{radius} which tells us when the John-Nirenberg radius is positive.

Some applications of the John-Nirenberg radius will be given. In Section 3, we will use the John-Nirenberg
radius to prove a well-known result: a positive harmonic function defined in a domain of a manifold with
a point removed is either a Green function, or smooth across the removed point.

In Section 4 and 5, we will apply John-Nirenberg radius to study a collapsing sequence of metrics in conformal geometry, i.e., we will show that, if $g_k=u_k^\frac{4}{n-2}g$ collapses, then there exists $c_k$ such that $c_ku_k$ converges to a positive function. Then, we show that the $\epsilon$-regularity in \cite{Li-Zhou} can be also deduced by employing John-Nirenberg radius. In Section 5, we will use the John-Nirenberg radius to prove that a sequence of metrics on a 4-dimensional manifold in a fixed conformal class with $\|K\|_{W^{1,2}}<C$ and fixed volume is compact in $C^{1,\alpha}$. The idea is, if the sequence blows up, then the neck domains can be considered as collapsing sequences. Then, by multiplying a suitable constant, one of the neck sequences converges to a complete flat manifolds with at least two ends collared topologically by $S^3\times \R$. Yet, this is impossible. Employing the same argument one can also give a new proof of the $C^{0,\alpha}$-compactness of a metric sequence, which is in a fixed conformal class and satisfies
$$
\vol(g_k)+\|K(g_k)\|_{L^p}<C,$$
where $p>\frac{n}{2}$. It is well-known that such a problem has been deeply studied by Chang-Yang \cite{Chang-Yang1, Chang-Yang2}, and solved by Gursky \cite{Gursky}.

\begin{center}
\begin{tikzpicture}

\draw (4,-4) arc (-66: 280: 1.3);
\draw (3.8,-4.25) arc (120: 450: .43 and .33);
\draw plot[smooth,tension=1.9]coordinates{(4,-4)(3.9,-4.1)(4.05,-4.2)};

\draw plot[smooth,tension=1.9]coordinates{(3.8,-4.25)(3.85,-4.15)(3.7,-4.1)};

\draw(4.2,-4.5) ellipse (.15 and 0.1);

\draw[blue!30](3.9,-4.15) ellipse (.15);

\draw (10,-4.5) arc (-50: 50: 2);
\draw (12,-1.4) arc (130: 230: 2);

\draw plot[smooth,tension=1.9]coordinates{(10.7,-3)(11,-3.1)(11.3,-3)};

\draw[dotted] plot[smooth,tension=1.9]coordinates{(10.7,-3)(11,-2.9)(11.3,-3)};

\draw[->,dashed]plot[smooth]coordinates{
(4.1,-4.1)(7,-3.75)(10.2,-3)};

\draw (11.9,-3) node {$S^{3}$};

\end{tikzpicture}

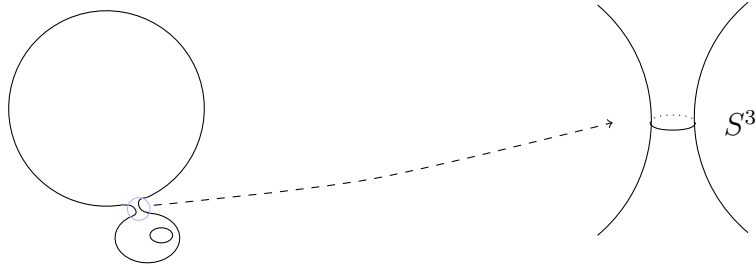
{\bf Fig 1.} After an appropriate rescaling, one of the
neck\\ sequences converges to a complete flat manifold,\\ which has at
least 2 ends collared topologically by
$S^3\times \R$.
\end{center}
\vspace{2ex}

In Section 6, we try to extend the definition of John-Nirenberg radius to the case of two dimensional manifolds. We will apply the John-Nirenberg radius to give a generalized H\'elein's Convergence Theorem. However, it is worthy to point out that Lemma \ref{radius.Lp.R} does not hold true for the case of two dimensional manifolds.

\section{John-Nirenberg radius}
First, we need to recall the following Moser's inequality on the ball $B^n$ for functions with mean value zero, which was established in \cite{Leckband}.

\begin{thm}\cite{Leckband}\label{M-T} Let $B_1$ be the unit ball of $\R^n$, and  $\alpha_n=n(\frac{\omega_{n-1}}{2})^{\frac{1}{n-1}}$, where $\omega_{n-1}$ is the measure of
unit sphere in $\mathbb{R}^n$. Then
$$
\sup_{u\in W^{1,n}(B_1),\,\, \int_Budx=0,\,\, \|\nabla u\|_{L^n(B_1)}\leq 1}\int_{B_1}e^{\alpha_n|u|^\frac{n}{n-1}}dx<+\infty.
$$
\end{thm}

From the Theorem above and the following inequality
$$
|u|=\frac{|u|}{\|\nabla_gu\|_{L^n}}\|\nabla_gu\|_{L^n}\leq
\frac{n-1}{n}\left(\frac{|u|}{\|\nabla u\|_{L^n}}\right)^\frac{n}{n-1}+
\frac{1}{n}\|\nabla u\|_{L^n}^n,
$$
we derive the following:

\begin{cor}\label{W-M-T} Let $B_1$ be the unit ball  of  $ \R^n$, and $
u\in W^{1,n}(B_1)$ and $\int_{B_{1}}u dx=0$.
Then
$$
\int_{B_1}e^{\beta_n|u|}dx<Ce^{\frac{\alpha_{n}}{n-1}\int_{B_1}|\nabla u|^ndx},
$$
where $\beta_n=\frac{n}{n-1}\alpha_n$.
\end{cor}

We say $u$ is essentially positive, if there exists $\epsilon>0$, such that
$u>\epsilon$ almost everywhere.
Given  an  essentially positive function $u\in W^{1,n}(\Omega)$, we define the John-Nirenberg radius as follows:
$$
\rho(x,u,\Omega,\epsilon_0)=\sup \left\{r:\int_{B_{r}(x)\cap\Omega}|\nabla \log u|^ndx<\epsilon^n_{0}\right\}.
$$
Later, $\rho(x,u,\Omega,\epsilon_0)$ will be used to study convergence of a sequence of positive functions. For example, we have the following:

\begin{lem}\label{example}
Let $\Omega$ be a domain of  $\R^n$, $u_{k}\in W^{1,n}(\Omega)$ be essentially positive. Let $\Sigma$ be a compact $(n-1)$-dimensional
submanifold (perhaps with boundary) embedded in $\Omega_1$, and $-\log c_k$ be the integral mean value of $\log u_k$ over $\Sigma$. Suppose
$\Omega_1\subset\subset\Omega$ and $\rho(x,u_k,\Omega,\epsilon_0)>a>0$ for any $x\in \Omega_1$. Then,
$c_{k}u_{k}$ and $\frac{1}{c_ku_{k}}$
are bounded in $L^{\frac{\beta_n}{\epsilon_0}}(\Omega_1)$. Moreover, $\log c_ku_k$
converges weakly in $W^{1,n}(\Omega_1)$.
\end{lem}

\proof Choose $a_1<\frac{1}{2}\min\{d(\Omega_1,\partial\Omega),a\}$, and define $\Omega_1'=\{x:d(x,\Omega_1)<a_1\}$.
By the assumptions, we have
$$\int_{\Omega_1'}|\nabla \log u_k|^n dx<C(\epsilon_{0},\Omega_{1}').$$
The Poincar\'e inequality tells us $\log c_{k}u_{k}$ is bounded in $W^{1,n}(\Omega_1')$,  where $-\log c_{k}$ is the integral mean value of $\log u_{k}$ on $\Sigma$.  Hence, we may assume that
$\log c_{k}u_{k}$ converges  in $L^{q}(\Omega_1')$ for any $q$.

Take $x_1$, $\cdots$, $x_m\in \overline\Omega_{1}$ such that $\{B_{a_1}(x_i): i=1, \cdots, m\}$ is an open cover of $\Omega_1$. Without loss of generality, we may assume
$\log u_k+\log c_k^i$ converges weakly in $W^{1,n}(B_{a_1}(x_i))$ and strongly in
$L^1(B_{a_1}(x_i))$. Here $-\log c_k^i$ is the mean value of $\log u_k$
over $B_{a_1}(x_i)$.
Since $$(\log u_k+\log c_k^i)-(\log u_k+\log c_k)$$ converges in $L^1(B_{a_1}(x_i))$,
we may assume $\log c_k^i-\log c_k$ converges.

By Corollary \ref{W-M-T}, we have
$$
\int_{B_{a_1}(x_i)}e^{\frac{\beta_n}{\epsilon_0}|\log u_k+\log c_k^i|}< C(\epsilon_{0},n),
$$
and hence
$$
\int_{B_{a_1}(x_i)}e^{\frac{\beta_n}{\epsilon_0}|\log u_k+\log c_k|}<C(\epsilon_{0},n).
$$
It turns out that both $\|c_ku_k\|_{L^\frac{\beta_n}{\epsilon_0}(\Omega_1)}$
and $\|\frac{1}{c_ku_k}\|_{L^\frac{\beta_n}{\epsilon_0}(\Omega_1)}$ are bounded.
\endproof

\begin{rem}
$-\log c_k$ in the above lemma can be chosen to be any constant which makes the Poincar\'e inequality
hold. For example, we can set $-\log c_k$ to be the mean value of $\log u_k$ over a subdomain of $\Omega_1$.
\end{rem}

We consider the operator
$$
L(u)=a^{ij}u_{ij}+b^iu_{i}+cu,
$$
where
\begin{equation}\label{elliptic.condition}
\|a^{ij}\|_{C^{0,\alpha}}+\|b^i\|_{C^{0,\alpha}}+\|c\|_{C^0}<
A_1,\s\s 0<A_2<a^{ij}<A_3.
\end{equation}
Later, we need to use the following:

\begin{cor}\label{case1}
Let $(\frac{1}{p}+\frac{2\epsilon_0}{\beta_n})<\frac{1}{2}+\frac{1}{n}$. Let $u_k\in W^{2,p}(B_2)$ be a sequence of positive functions, each of which solves the equation
$Lu_k=f_ku_k$ where $\|f_k\|_{L^{p}(B_2)}<\Lambda$.
If
$$
\|u_k\|_{L^\frac{\beta_n}{\epsilon_0}(B_2)}+\|\frac{1}{u_k}\|_{L^\frac{\beta_n}{\epsilon_0}(B_2)}<\Lambda_2,
$$
 then, after passing to a subsequence,
$u_k$ converges weakly in $W^{2,q}(B_1)$ and
$\log u_k$ converges weakly in $W^{2,q'}(B_1)$ for any
$$
q\in\left(\frac{1}{\frac{n+2}{2n}-\frac{\epsilon_{0}}{\beta_{n}}},\, \frac{p}{1+\frac{\epsilon_{0}}{\beta_{n}}p}\right)\cap (1,\,n)\s\s
\mbox{and}\s\s q'\in \left(1,\, \frac{1}{2(\frac{\epsilon_0}{\beta_n}+\frac{n-q}{nq})}\right)
\cap (1,\, p).
$$
\end{cor}

\proof
Since $\frac{p}{1+\frac{\epsilon_{0}}{\beta_{n}}p}>q$, we have $\frac{pq}{p-q}<\frac{\beta_n}{\epsilon_0}$. Noting
$$
\int_{B_2}|f_ku_k|^qdx\leq\left(\int_{B_2}|f_k|^pdx\right)^\frac{q}{p}
\left(\int_{B_2}|u_k|^\frac{pq}{p-q}dx\right)^\frac{p-q}{p},
$$
by the standard elliptic theory we get the estimate of $\|u_k\|_{W^{2,q}(B_\frac{3}{2})}$. Then, it follows
$$\|\nabla u_k\|_{L^\frac{nq}{n-q}(B_\frac{3}{2})}<C.$$
It is easy to check that $$2q'<\frac{nq}{n-q} \s\s\s\mbox{and}\s\s\s \frac{2nqq'}{nq-2q'(n-q)}<\frac{\beta_n}{\epsilon_0}.$$
By H\"{o}lder inequality, we have
$$
\int_{B_\frac{3}{2}}|\nabla\log u_k|^{2q'}dx<
\left(\int_{B_\frac{3}{2}}|\nabla u_k|^\frac{nq}{n-q}dx
\right)^\frac{2q'(n-q)}{nq}\left(\int_{B_\frac{3}{2}}
\left(\frac{1}{u_k}\right)^\frac{2nqq'}{nq-2q'(n-q)}dx\right)^\frac{nq-2q'(n-q)}{nq}<C_{1}.
$$
Define an operator $L'=L-c$. Obviously, $\log u_{k}$ satisfies the following equation
$$
L'(\log u_k)=-a^{ij}(\log u_k)_i(\log u_k)_{j}+f_k-c.
$$
By $L^{p}$ estimate, we know $\|\log u_k\|_{W^{2,q'}}<C_{2}$.
\endproof

\begin{rem}
In Corollary \eqref{case1}, in order to guarantee that $$\frac{p}{1+\frac{\epsilon_{0}}{\beta_{n}}p}> \frac{1}{\frac{n+2}{2n}-\frac{\epsilon_{0}}{\beta_{n}}},$$
we only need to choose $p$ such that $\frac{1}{p}+\frac{2\epsilon_{0}}{\beta_{n}} < \frac{1}{2}+ \frac{1}{n}$.
Hence, it follows that $$\left(\frac{1}{\frac{n+2}{2n}-\frac{\epsilon_{0}}{\beta_{n}}},\,\,\frac{p}{1+\frac{\epsilon_{0}}{\beta_{n}}p}\right)\neq \emptyset \s\s\s \mbox{and}\s\s\s
\left(1,\,\, \frac{1}{2(\frac{\epsilon_{0}}{\beta_{n}}+\frac{n-q}{nq})}\right)\neq \emptyset.$$
\end{rem}

The following theorem is the key point of this paper:
\begin{thm}\label{radius}
Let $p>\frac{n}{2}$,
and $\frac{1}{p}+\frac{2\epsilon_{0}}{\beta_n}<\frac{2}{n}$.
Let $u_k\in W^{2,p}(B_3)$ be a smooth positive function which solves the equation
$$
Lu_k=f_ku_k.
$$
If for any $x_k\rightarrow x_0\in \overline{B_2}$
and
$r_k<2\rho(x_k,u_k,B_3,\epsilon_0)$ with $r_k\rightarrow 0$, a subsequence of $r_k^2f_k(r_kx+x_k)$ is bounded in
$L^p(B_\frac{1}{4})$, and converges to 0 in the sense of
distribution on $B_\frac{1}{4}$,
then there exists $a>0$, such that
$$
\rho(x,u_k,B_3,\epsilon_0)>a,\s\s \forall x\in B_1.
$$
\end{thm}

\proof We argue by contradiction. Assume the conclusion is not true. Then we can find $x_k\in B_1$, s.t.
$\rho(x_k,u_k,B_3,\epsilon_0)\rightarrow 0$.
For simplicity, we denote $\rho(x,u_k,B_3,\epsilon_0)$ by $\rho_k(x)$.

Put $y_k\in B_2$, such that
$$
\frac{\rho_k(y_k)}{2-|y_k|}=\inf_{x\in B_2}\frac{\rho_k(x)}{2-|x|}:=\lambda_k.
$$
Noting that $$\lambda_k\leq \frac{\rho_k(x_k)}{2-|x_k|}\leq  \rho_{k}(x_k)\rightarrow 0,$$
we have $\rho_k(y_k)\rightarrow 0$, and hence for any fixed $R$
$$
\frac{\rho_k(y_k)}{2-|y_k|}\rightarrow 0, \s\s\s
B_{R\rho_k(y_k)}(y_k)\subset  B_{2-|y_{k}|}(y_{k})\subset B_2,
$$
when $k$ is sufficiently large.
Then, for any $y\in B_{R\rho_k(y_k)}(y_k)$ we have
\begin{eqnarray*}
\frac{\rho_k(y)}{\rho_k(y_k)}&\geq& \frac{2-|y|}{2-|y_k|}\geq \frac{2-|y_k|-|y-y_k|}{2-|y_k|}\\
&\geq&1-\frac{R\rho_k(y_k)}{2-|y_k|}\\
&=&1-R\lambda_k.
\end{eqnarray*}
Hence, as $k$ is large enough, there holds
 $$\frac{\rho_k(y)}{\rho_k(y_k)}>\frac{1}{2}.$$

Assume $y_k\rightarrow y_0$.
Let $v_k(x)=u_k(y_k+r_kx)$, where $r_k=\rho_k(y_k)$. Then, there holds
$$\rho(x,v_k,B_R,x)\geq 1/2\s\s \mbox{on}\s\s B_\frac{R}{2}$$
and
$$
\int_{B_1}|\nabla\log v_k|^n=\epsilon_0^n.
$$
Moreover, $v_k$ satisfies the following equation:
$$
a^{ij}(y_k+r_kx) (v_k)_{ij}=-r_kb^i(y_k+r_kx)(v_k)_i-cr_k^2v_k+
r_k^2f_k(y_{k}+r_kx)v_k.
$$
By Lemma \eqref{example}, we can find $c_k$, such that
$$
\|c_kv_k\|_{L^\frac{\beta_n}{\epsilon_0}( B_{\frac{R}{2}})}+
\|\frac{1}{c_kv_k}
\|_{L^\frac{\beta_n}{\epsilon_0}( B_\frac{R}{2})}<C(R).
$$

Noting that $\{r_k^2f_k(y_{k}+r_{k}x_{0}'+r_{k}x)\}$ is bounded in $L^p(B_\frac{1}{4})$ for any $x_0'\in B_\frac{R}{2}$, by a covering argument we can see that the sequence $\{r_{k}^{2}f_k(y_k+r_{k}x)\}$ is bounded in $L^p(B_{R/2})$ for any $R$. By the same arguments, we also know that in the sense of distribution on $\R^n$
$$r_{k}^{2}f_k(y_{k}+r_{k}x)\rightarrow 0.$$

By the assumptions, we have
$$
2(\frac{\epsilon_0}{\beta_n}+\frac{n-q}{nq})<\frac{2}{n}
$$
when $q=\frac{p}{1+\frac{\epsilon_0}{\beta_n}p}$. Noting that $\frac{2}{n}<\frac{1}{2}+\frac{1}{n}$,
by Corollary \ref{case1}, we can find $q$ and $q'>\frac{n}{2}$, such that $c_kv_k$ converges to a function $v$ weakly in $W^{2,q}(B_R)$,
and $\log c_kv_k$  converges weakly  to $\log v$  in $W^{2,q'}(B_R)$. Then $\log c_kv_k$ converges in $C^0(R)$, which implies that $v>0$.
Moreover, $a^{ij}(y_0)v_{ij}=0$, then $v$ is a positive constant. However, by the Sobolev embedding theorem,
$|\nabla \log c_kv_k|$ converges in $L^n$.
Then, it follows
$$
\int_{B_1}|\nabla \log v|^n=\epsilon_0^n,
$$
which is impossible since $\log v$ is a constant. Thus we complete the proof.
\endproof

\begin{cor}\label{radius.2}
Let $p$, $\epsilon_0$ be as in Theorem \ref{radius}.
Let $u\in W^{2,p}(B_3)$ be a smooth positive function which solves the equation
$$
Lu=fu.
$$
Then there exist positive numbers  $\epsilon$ and  $a$ which only depend on
$A_1$, $A_2$, $A_3$, $p$ and $\epsilon_0$ such that, if
\begin{equation}\label{harmonic.2}
\sup_{x\in B_2,r<2\rho(x,u,B_3,\epsilon_0)}r^{2p-n}\int_{B_r(x)}|f|^{p}<\epsilon,
\end{equation}
then
$$
\rho(x,u,B_3,\epsilon_0)>a,\s\s \forall x\in B_1.
$$
\end{cor}

\proof We argue by contradiction. Assume the above conclusion is not true. Then there exists a sequence of $u_k$  satisfying $$Lu_k=f_k,$$ such
that
$$
\inf_{B_1}\rho(x,u_k,B_3,\epsilon_0)\rightarrow 0
$$
and
$$
\sup_{x\in B_2,r<2\rho(x,u_{k},B_3,\epsilon_0)}r^{2p-n}\int_{B_r(x)}|f_k|^{p}\rightarrow 0.
$$
It is easy to check from the above that $r_k^2 f_k(x_k+r_{k}x)$ converges to 0 in $L^p(B_{\frac{1}{4}})$. Thus we get the desired conclusion from Theorem \ref{radius}.
\endproof

\section{Positive harmonic function with isolated singularity}
In this section, we will use the so-called John-Nirenberg radius or the John-Nirenberg inequality to study the positive  harmonic   functions with singularity on a manifold. We will prove the following:

\begin{lem}\label{harmonic}
Let $g=dr^2+g(r,\theta)d
\mathbb{S}^{n-1}$ be a smooth metric over $B_1\subset\R^n$, where $g(r,\theta)=r^2(1+o(1))$.
Assume $u$ is a positive harmonic function on $B_{1}\setminus\{0\}$. Then  $u\in W^{1,q}$ for any  $q\in (1,\frac{n}{n-1})$
and satisfies the weak equation
$$
-\Delta_g u=c\delta_0,\s \s c\geq 0.
$$
\end{lem}

\proof
Let $$c=-\int_{\partial B_r}\frac{\partial u}{\partial r}dS_r.$$
First, we prove that
$$
\frac{u(rx)}{r^{2-n}}\rightarrow \frac{c}{(n-2)\omega_{n-1}}\s\s and\s\s \frac{\nabla u}{r^{1-n}}(rx)\rightarrow -\frac{c}{
\omega_{n-1}}\frac{\partial}{\partial r}
$$
uniformly on $S^{n-1}$.

Assume this is not true. Then we can find $x_k\in S^{n-1}\subset\R^n$ and
$r_k\rightarrow 0$, such that
$$
\left|\frac{u(r_kx_k)}{r_k^{2-n}}-\frac{c}{(n-2)\omega_{n-1}}\right|>\epsilon\s\s\s \mbox{or}\s\s\s  \left|\frac{\nabla u}{r_{k}^{1-n}}(r_kx_k)+\frac{c}
{\omega_{n-1}}\frac{\partial}{\partial r}\right|>\epsilon.
$$
Let  $v_k=u(r_kx)$ and choose $c_{k}$ such that $$\int_{\partial B_{1}}\log c_kv_k d\CS^{n-1}=0.$$
By the results in the above section, for any $r>0$ we can find $a(r)>0$ such that, for any $x\in B_\frac{1}{r}\setminus B_r$,
$$\rho(x,v_k,B_{\frac{2}{r}}\setminus B_\frac{r}{2}, \epsilon_{0})>a(r),$$
and hence both $c_kv_k$ and $\frac{1}{c_kv_k}$ are bounded in $L^{\frac{\beta_n}{\epsilon_0}}(B_{\frac{1}{r}}\setminus B_r)$.\medskip

Since $c_kv_k$ is harmonic, after passing to a subsequence, $c_kv_k$ converges
in $C^\infty_{loc}(\R^n)$
to a  function $v$ which is positive and harmonic on $\R^n\setminus
\{0\}$. It is well-known that
$$v=ar^{2-n}+b,$$
where $a$ and $b$ are nonnegative real numbers with $a^{2}+b^{2} > 0$ (cf Corollary 3.14 \cite{Sheldon-Paul-Wade}).

Now, we need to discuss the following two cases.

Case 1: $c\neq 0$. In this case, from
\begin{equation}\label{candck}
\frac{c_k}{r_{k}^{n-2}}\int_{\partial B_{r_k}}\frac{\partial
u}{\partial r}dS_r\rightarrow\int_{\partial B_1}\frac{\partial v}
{\partial r}d\CS^{n-1}=a(2-n)\omega_{n-1},
\end{equation}
it follows that
\begin{equation}\label{cneq0}
\frac{c_k}{r_{k}^{n-2}}\rightarrow a\frac{(n-2)\omega_{n-1}}{c}.
\end{equation}
Then we have
$$
\frac{u_k(r_kx)}{r_k^{2-n}}=\frac{c_kv_k(x)}{c_kr_k^{2-n}}\rightarrow\frac{(a+b)c}{a(n-2)\omega_{n-1}},
$$
and
$$
r_k^{n-1}\nabla u(r_kx)\rightarrow-\frac{c}{\omega_{n-1}}\frac{\partial}{\partial r}.
$$
To get a contradiction, we need to prove $b=0$. Let $G$ be the Green function which satisfies
$-\Delta_g G=\delta_0$ and $G|_{\partial B _\delta}=0$.
We have
$$
\lim_{r\rightarrow 0}r^{n-2}G=\lim_{r\rightarrow 0}
\frac{r^{n-1}}{2-n}\frac{\partial G}{\partial r}=\lambda\neq 0.
$$
Then
\begin{eqnarray}\label{b=0}
c_{k}\int_{\partial B_t}(u\frac{\partial G}{\partial r}-\frac{\partial  u}
{\partial r}G)dS_g&=&c_{k}\int_{\partial B_{r_k}}(u\frac{\partial G}{\partial r}-\frac{\partial u}
{\partial r}G)dS_g\nonumber\\
&=&\int_{\partial B_1}(\lambda c_{k}v_{k}(x)(2-n)r_{k}^{1-n}(1+o(1))-\lambda \frac{\partial c_{k}v_{k}}{\partial r}r_{k}^{1-n})\nonumber\\
&&\hspace{5ex}\times r_{k}^{n-1}(1+o(1))dS\nonumber\\
&\rightarrow&\int_{\partial B_1}\left(v(1)(2-n)\lambda
-v'(1)\lambda \right)dS\\
&=&-b(n-2)\omega_{n-1}\lambda .\nonumber
\end{eqnarray}
Obviously \eqref{cneq0} implies $c_{k}\rightarrow 0$, hence from the above equality \eqref{b=0} we derive that $b=0$.\medskip

Case 2: $c=0$. If $c_kr_k^{2-n}\rightarrow+\infty$, it is easy to check that
$$r_k^{n-2}u(r_kx)=\frac{c_kv_k(x)}{c_kr_k^{2-n}}\rightarrow 0\s\s\s \mbox{and}\s\s\s r_k^{n-1}\nabla u(r_kx)\rightarrow 0.$$
On the other hand, if $c_kr_k^{2-n}<C$, it follows that $c_k\rightarrow 0$.
From \eqref{candck} we have $a=0$. From \eqref{b=0} we can see that $b=0$. Thus, we get a contradiction.

Therefore, we conclude that $u\in W^{1,q}(B)$ for any $q\in(1,\frac{n}{n-1})$. Given a smooth function $\varphi$ whose support set is contained in $B_1$, we have
\begin{eqnarray*}
\int_{B_1}\nabla\varphi\nabla u dV_g&=&
\lim_{r\rightarrow 0}\int_{B_1\setminus B_r}\nabla\varphi\nabla u dV_g\\
&=&-\lim_{r\rightarrow 0}\int_{B_1\setminus B_r}\Delta u\varphi dV_g
+\lim_{r\rightarrow 0}\int_{\partial B_r}
\frac{\partial u}{\partial \nu}\varphi dV_g\\
&=&
c\varphi(0).
\end{eqnarray*}
Thus, we get $$-\Delta_g u=c\delta_0.$$
Thus we complete the proof of this lemma.
\endproof

\begin{cor}\label{singularity}
Let $(M,g)$ be a closed manifold with constant scalar curvature
$R(g)$. Suppose $p_1$, $\cdots$, $p_m\in M$,
and $g'$ is a metric on $M\setminus\{p_1,\cdots,p_m\}$, which
is conformal to $g$. If $R(g')=0$, then $(M,g')$ is complete near
$p_i$ or $g'$ is smooth across
$p_i$.
\end{cor}

\proof
We can find a metric $g_0$ which is conformal to $g$, such that
$R(g_0)=0$ in a neighborhood of  $p_i$. Let $g'=u^\frac{4}{n-2}g_0$.
Then $u$ is harmonic in a neighborhood of $p_i$. Thus, either $u$
can be extended smoothly to $p_i$, or $u\sim c_ir^{2-n}$
for a positive $c_i$, which implies that $g'$ is complete
near $p_i$.
\endproof

\section{A  collapsing sequence  with bounded $\|R\|_{L^p}$}
In the previous paper \cite{Li-Zhou},  the authors  use the $\epsilon$-regularity to study the bubble tree
convergence of a metric sequence in a fixed conformal class with bounded  volume and $L^p(M)$-norm of scalar curvature.
Then, it has been shown that the John-Nirenberg radius is bounded below by a positive constant when the
volume converges to 0.
In this section, we will show that the $\epsilon$-regularity is also a corollary of John-Nirenberg inequality, which was deduced
directly from $L^p$ estimate in \cite{Li-Zhou}.

First, we show the positivity of the John-Nirenberg radius for a
collapsing sequence.

\begin{lem}\label{radius.Lp.R}
Let $n>3$ and  $\{\hat{g}_k\}$ be a sequence of  metrics over $B_2\subset\R^n$ which converges to $g$.
Let $g_k=u_k^\frac{4}{n-2}\hat{g}_k$. Suppose that $\vol(B_2,g_k)\rightarrow 0$ and
$\int_{B_2}|R(g_k)|^pd\mu_{g_k}<\Lambda$, where
$p>\frac{n}{2}$. Then for any sufficiently small $\epsilon_0$, there exists $a_0>0$, such that
$$
\rho(x,u_k,B_2,\epsilon_0)>a_0,\s \forall\ x\in B_1.
$$
\end{lem}

\proof
Let $$f_k=R(g_k)u_k^\frac{4}{n-2}.$$
Given  $x_k\rightarrow x_0$, $r_k\rightarrow 0$, such that
$$
r_k<2\rho(x_k,u_k,B_2,\epsilon_{0}).
$$
We let $v_k(x)=r_k^\frac{n-2}{2}u_k(x_k+r_kx)$. Obviously, $\rho(0,v_k,B_2,\epsilon_{0})\geq 1/2$ implies $$\int_{B_\frac{1}{3}}|\nabla\log v_k|^ndx<\epsilon_0.$$
Then, $$\rho(y,v_{k},B_2(0),\epsilon_{0})>\frac{1}{6}, \s\s \forall \ y\in B_{\frac{1}{3}}(0).$$
By Lemma \ref{example},
$$
\|c_kv_k\|_{L^\frac{\beta_n}{\epsilon_0}(B_\frac{1}{3})}+\|\frac{1}{c_kv_k}\|_{L^\frac{\beta_n}{\epsilon_0}(B_\frac{1}{3})}\leq C.
$$
Since $\int |v_k|^\frac{2n}{n-2}\rightarrow 0$, we get
$c_k\rightarrow+\infty$, then $\|v_k\|_{L^\frac{\beta_n}{\epsilon_0}(B_\frac{1}{3})}<C$. Fix a $p'\in(\frac{n}{2},p)$, we can
choose $\epsilon_0$ to be sufficiently small such that
such that
$$
\|R(\hat{g}_k)(x_k+r_kx)v_k^\frac{4}{n-2}\|_{L^{p'}(B_\frac{1}{3})}<
C(p',p,\Lambda),\s and\s
\frac{p'}{1+
\frac{\epsilon_0}{\beta_n}p'}>\frac{n}{2}.
$$
Applying Corollary \ref{case1} to $c_kv_k$, we get
$c_kv_k$ converges to a positive function $\phi$ in $W^{2,q}(B_\frac{1}{4})$
for some $q>\frac{n}{2}$. Since $\vol(B_{2},g_{k})\rightarrow 0$,  $c_{k}\rightarrow \infty$, and
$v_{k}\rightarrow 0$ in $B_{\frac{1}{4}}$  uniformly, we get
\begin{eqnarray*}
\int_{B_{\frac{1}{4}}(0)}r_{k}^{2p}|R(\hat g_{k})|^{p}u_{k}(x_{k}+r_{k}x)^{\frac{4p}{n-2}}dx
&=&\int_{B_{ \frac{1}{4}}(0)}|R(\hat g_{k})|^{p}v_{k}^{\frac{4p}{n-2}}dx\\
&\leq&||v_{k}||_{L^{\infty}}^{\frac{4p-2n}{n-2}}\int_{B_{ \frac{1}{4}}(0)}|R(\hat g_{k})|^{p}v_{k}^{\frac{2n}{n-2}}dx\rightarrow 0
\end{eqnarray*}
Applying Theorem \ref{radius}, we obtain the required result and complete the proof.
\endproof

Next, we  prove the  $\epsilon$-regularity.

\begin{lem}
Let $B_r(x_0)\subset M$. We assume $\int_{B_r(x_0)}|R(g')|^pd\mu_{g'}<\Lambda$, where
$p>\frac{n}{2}$. Then, there exists $\epsilon_0'$ which depends only on $M$, $r$ and
$\Lambda$, such that if $\vol(B_r(x_0),g')<\epsilon_0'$, then
$$
\|u\|_{W^{2,p}(B_\frac{r}{2})(x_0)}<C\vol(B_r(x_0),g')^\frac{n-2}{2n}.
$$
\end{lem}
\proof Assume the result is not true. Then we can find $x_k\rightarrow x_0'$,
$g_k=u_k^\frac{4}{n-2}g$, such that $\vol(B_r(x_k),g_k)\rightarrow 0$,
$\|R(g_k)\|_{L^p(B_r(x_k))}<\Lambda$ and
$$
\|u_k\|_{W^{2,p}(B_\frac{r}{2}(x_0))}>k\vol(B_r(x_0), g_{k})^\frac{n-2}{2n}.
$$
$g_k'=g|_{B_r(x_k)}$ can be regarded as a metric over $B_r\subset\R^n$ which converges smoothly. It follows that
$$
-\Delta_{g_k'}u_k=(-c(n)R(g_k')+R(g_k)u_k^\frac{4}{n-2})u_k.
$$

By the above lemma, $\rho(x,u_k,B_r,\epsilon_0)>a>0$ for any $ x\in B_{\frac{7r}{8}}$.
Choose $\epsilon_0$ to be sufficiently small.
By Corollary \ref{case1},
$c_ku_k$ converges in $W^{2,q}(B_{\frac{3}{4}r})$ to a positive function for some $q>\frac{n}{2}$. Then,
$\|c_ku_k\|_{L^\infty(B_{\frac{3}{4}r})}$ is bounded above. Since $\int_{B_r(x_0)}u_k^\frac{2n}{n-2}
\rightarrow 0$, we get $c_k\rightarrow+\infty$. Then
 $u_{k}\rightarrow 0$ in $B_{\frac{r}{2}}(x_{0})$ uniformly and
$$
\|c_ku_k\|_{W^{2,p}(B_\frac{r}{2}(x_0))}
\geq k\vol(B_{\frac{3}{4}r}(x_0), c_k^\frac{4}{n-2}g_{k})^\frac{n-2}{2n}
\rightarrow+\infty.
$$

On the other hand, since
$$
\int_{B_{\frac{3}{4}r}}(|R(g_k)|u_k^\frac{4}{n-2}c_ku_k)^p\leq
\int_{B_{\frac{3}{4}r}}|R(g_k)|^pu_k^\frac{2n}{n-2}
\|(c_ku_k)^{\frac{p(n+2)-2n}{n-2}}\|_{L^\infty(B_{\frac{3}{4}r})}c_k^{p-\frac{p(n+2)-2n}{n-2}}
\rightarrow 0,
$$
we derive
\begin{eqnarray*}
\|c_ku_k\|_{W^{2,p}(B_\frac{r}{2})}&\leq& C(\|R(g_k)u_{k}^{\frac{4}{n-2}}c_{k}u_{k}\|_{L^p(B_{\frac{3}{4}r})}+\|c_ku_k\|_{L^p(B_{\frac{3}{4}r})})\\
&\leq& C(\|R(g_k)\|_{L^p(B_r,g_k)}+1)\\
&<&C.
\end{eqnarray*}
We get a contradiction and finish the proof.
\endproof

\section{4-manifold in a conformal class with $\|K\|_{W^{1,2}}<\Lambda$}
In this section, we let $\dim M=4$, $u_{k}\in W^{3,2}(M,g)$ and $g_k=u_k^2g$.
Assume that for every $k$ there hold
\begin{equation}\label{4dim}
\vol(M,g_k)=1\s\s \mbox{and} \s\s \int(|\nabla_{g_k}K(g_k)|^2+K^2(g_k))d\mu_{g_{k}}<\Lambda,
\end{equation}
where  $K(g_k)$ denotes the sectional curvature of $g_k$. We intend to study the convergence behavior of $u_k$.

First of all, we try to show that the John-Nirenberg inequality will imply the $L^p$-estimate of curvature. We want to prove the result under the assumption that
$$\|R(g_k)\|_{W^{1,2}(M,g_k)}<\Lambda.$$

\begin{lem}\label{4dim.convergence}
Let $g=g_{ij}dx^i\otimes dx^j$ be a smooth metric on $B_3\subset\R^n$ with
$\|g_{ij}\|_{C^{2,\alpha}(B_3)}<\gamma_1$. Suppose that $g'=u^2g$ satisfies $\vol(B_3,g')<\gamma_2$ and
\begin{equation}\label{R-4dim}
\int_{B_3}(|\nabla_{g'}R(g')|^2+|R(g')|^2)
dV_{g'}<\Lambda.
\end{equation}
Then, for any $p<4$, there exists $\hat{\epsilon}_0=\hat\epsilon_0(p)$ such that, if $\epsilon_0<\hat{\epsilon}_0$  and
$\rho(x,u,B_3,\epsilon_{0})\geq a>0$, there holds true
$$
\int_{B_1}|R(g')u^2|^p<C(p,\gamma_1,\gamma_2,\hat{\epsilon}_0,\Lambda,a).
$$
\end{lem}

\proof For any $q\in(\frac{4}{3},2)$,
we have
\begin{eqnarray*}
\int_{B_1}|\nabla_g(Ru^2)|^qdV_g
&\leq& C(q)\left(\int_{ B_1}(|\nabla_gR|u)^qu^qdV_g+
\int_{B_1}(|R|u^2)^q|\nabla_g\log u|^qdV_g\right)\\
&\leq&  C(q,\Lambda)\left(\int_{B_1}u^\frac{2q}{2-q}dV_g\right)^\frac{2-q}{2}\\
&&+
C(q)\left(\int_{B_1}(|R|u^2)^\frac{4q}{4-q}\right)^\frac{4-q}{4}\left(\int_{B_1}|\nabla_g\log u|^4dV_g\right)^{^\frac{q}{4}}.
\end{eqnarray*}

Choose $\hat{\epsilon}_0$, such that $\frac{\beta_n}{\hat\epsilon_0}>\frac{2q}{2-q}$.
Let $p=\frac{4q}{4-q}$,  and   $-\log c$ be the mean value of $\log u$ over $B_1$.
By Corollary \ref{W-M-T}, we can find  $C=C(\epsilon_{0},q,\gamma_{1},a)$, such that
both   $\|cu\|_{L^\frac{2q}{2-q}(B_{1})}$ and $\|\frac{1}{cu}\|_{L^\frac{2q}{2-q}(B_{1})}$ are bounded above by $C$. Then
$$
|B_1|^2\leq \int_{B_1}(cu)^4dx
\int_{B_1}(cu)^{-4}dx<C(\epsilon_{0},q,\gamma_{1},a)c^4\int_{B_1}u^4,
$$
which yields that $c$ is bounded below by a positive constant  $C=C(\epsilon_{0},q,\gamma_{1},\gamma_{2},a)$. Then
\begin{eqnarray*}
\int_{B_1}|\nabla(Ru^2)|^q  dV_{g}
&\leq& C(q,\epsilon_{0},\gamma_{1},\gamma_{2},\Lambda,a)+   C(q,\gamma_1,a)\epsilon_0^{q}\left(\int_{B_1}(|R|u^2)^\frac{4q}{4-q}dV_g\right)^\frac{4-q}{4}.
\end{eqnarray*}
Let  $\epsilon=C(q,\gamma_1,a)^{\frac{1}{q}}\epsilon_0$. We get
$$
\|\nabla(Ru^2)\|_{L^q(B_1,g)}<  C(q,\epsilon_{0},\gamma_{1},\gamma_{2},\Lambda,a) +\epsilon\|Ru^2\|_{L^\frac{4q}{4-q}(B_1,g)}.
$$
By Sobolev inequality,
$$
\|Ru^2\|_{L^\frac{4q}{4-q}(B_1,g)}\leq  C(q,\gamma_1)(\|\nabla(Ru^2)\|_{L^q(B_1,g)}+
\|Ru^2\|_{L^q(B_1,g)}).
$$
Put $\epsilon C(q,\gamma_1)<\frac{1}{2}$, we get
\begin{eqnarray*}
\|Ru^2\|_{L^\frac{4q}{4-q}(B_1)}
&\leq&  C(q,\epsilon_{0},\gamma_{1},\gamma_{2},\Lambda,a).
\end{eqnarray*}
\endproof

Next, we show that $\|R\|_{L^2}$ small implies the
boundness of John-Nirenberg radius.
\begin{lem}\label{epsilon}
Let $g$,  $u$ and $g'$ be as in the above lemma. Then, there exist $\tau>0$ and $a>0$ such that, if $\int_{B_3}R^2(g')dV_{g'}<\tau$, then
$$
\inf_{x\in B_1}\rho(x,u,B_3,\epsilon_{0})>a,\s\s \forall \ x\in B_{1}.
$$
\end{lem}

\proof
We prove it by contradiction.
Assume there exists $g_{k}=u_{k}^{2}g$  such that
$$
\inf_{x\in B_1}\rho(x, u_{k},  B_3, \epsilon_{0})\rightarrow 0.
$$
Given $y_{k}\rightarrow y_{0}\in  \overline{B_2}$, $r_{k}<2\rho(y_{k},u_{k},B_{3},\epsilon_{0})$, we set $v_{k}(x)=r_{k}u_{k}(y_{k}+r_{k}x)$ and
$$
\hat{g}_k=v_k^2g_{ij}(y_k+r_kx)dx^i\otimes dx^j.
$$
Then, it is easy to see that $\rho(0,v_{k}, B_\frac{1}{2},\epsilon_{0})> \frac{1}{2}$ and $\rho(x,v_{k}, B_\frac{1}{2},\epsilon_{0})> \frac{1}{4}, \forall x\in B_{\frac{1}{4}}$.

Obviously,
$$
\|R(\hat g_{k})\|_{W^{1,2}(B_3,\hat{g}_k)}=\|R(  g_{k})\|_{W^{1,2}(B_{3r_k}(y_k),g_k)}.
$$
By Lemma \ref{4dim.convergence}, for some $p\in (2,\, 4)$ there holds
 $$\int_{B_{\frac{1}{4}}}|R(\hat g_{k})v_{k}^{2}|^{p}<C.$$
Since
$$
\int_{B_{\frac{1}{4}}}|r_{k}^{2}R(g_{k})(r_{k}x+y_{k})u_{k}^{2}(r_{k}x+y_{k})|^{p}dx\leq C\int_{B_{\frac{1}{4}}}|R(\hat g_{k})v_{k}^{2}|^{p}dx,
$$
and
$$
\int_{B_{\frac{1}{4}}}|r_{k}^{2}R(r_{k})(r_{k}x+y_k)u_{k}^{2}(y_{k}+r_{k}x)|^{2}=\int_{B_{\frac{1}{4}}}|R(\hat g_{k})v_{k}^{2}|^{2}=\int_{B_{\frac{1}{4}{r_{k}}}(y_{k})}|R(g_{k})|^{2}u_{k}^{4}\rightarrow 0.
$$
From Lemma \ref{radius}, it follows that $\rho(x,u_{k},B_{3},\epsilon_{0})> a, \forall x\in B_{1}$. Then, we get a contradiction.
\endproof\\

For convenience, given a subset $A\subset \mathbb{S}^{n-1}$, we set
$$
A_r=\bigcup_{t\in(0,r]}tA,\s\s\s\s C(A,r)=\bigcup_{t\in [\frac{r}{2},r]} tA.
$$
We need to establish the following lemma:

\begin{lem}\label{sing1}
Let  $g$ be a smooth metric over $B_1\subset\R^4$ and $g'=u^2g$, where $u\in W^{3,2}(B_1)$ is a positive function. Assume $g=dr^2+g(r,\theta)d
\mathbb{S}^{3}$ with  $g(r,\theta)=r^2(1+o(1))$. If
$$
\vol(B_1,g')+\int_{B_1}(|K(g')|^2+|\nabla_{g'}K(g')|^2)dV_{g'}<+\infty,
$$
then, when $r$ is small enough, there holds
$$
\vol(C(A,r/2),g')<\frac{1}{2^3} \vol(C(A,r),g').
$$
\end{lem}

\proof
We claim that: there exists $r_0$, such that if $r<r_0$, then
\begin{equation}\label{3circles1}
\vol(C(A,r/4),g')<\frac{1}{2^3}\vol(C(A,r/2),g')<\frac{1}{2^6}\vol(C(A,r),g'),
\end{equation}
or
\begin{equation}\label{3circles2}
\vol(C(A,r),g')<\frac{1}{2^3} \vol(C(A,r/2),g')<\frac{1}{2^6}  \vol(C(A,r/4),g').
\end{equation}

Assume there exists $r_k\rightarrow 0$, such that none of the above holds.
Put $u_{k}(x)=r_{k}u(r_{k}x)$ and $g_{k}=u_{k}^{2}g(r_kx)$. For any fixed
$R$, we have
$$
\int_{B_R}(|K( g_{k})|^{2}+|\nabla_{g_k}K(g_k)|^2)dV_{ g_{k}}=\int_{B_{Rr_{k}}}(|K(g')|^{2}+|\nabla_{g'}K(g')|^2)dV_{g'}\rightarrow 0.
$$
Then by Lemma \ref{4dim.convergence}-\ref{epsilon}, Lemma \ref{case1} and Lemma \ref{example}, we can find $\tilde c_{k}$ such that
$\tilde c_{k}u_{k}$ converges to a positive function $\varphi$. Let
$\tilde{g}_k=\tilde c_k^2g_k$ and $\tilde{g}=\varphi^2g_{\R^4}$.

Since $\vol(g_{k},B_{Rr_k}\setminus\{0\})\rightarrow 0$,
we have $\tilde c_{k}\rightarrow \infty$. Then it is easy to check that
$$
\int_{B_\frac{1}{r}\setminus B_r}(|\nabla_{\tilde{g_k}}K(\tilde{g}_k)|^2+
|K(\tilde{g}_k)|^2)dV_{\tilde{g}_k}\rightarrow 0.
$$
By Lemma \ref{epsilon} again, $\tilde c_ku_k$ converges weakly in $W^{3,2}_{loc}(\R^4\setminus\{0\})$ and $K(\varphi)=0$, thus $\varphi$ is a positive harmonic function.

Theorem 9.8 in \cite{Sheldon-Paul-Wade} tells us that $\varphi$ can be written as $\varphi = ar^{-2}+b$. Since $K(\varphi)=0$, we get $a=0$ or $b=0$.
When $a=0$ and $b\neq 0$, we have
$$
\frac{\vol(C(A,1),\tilde{g})}{\vol(C(A,1/2),\tilde{g})}=
\frac{\vol(C(A,1/2),\tilde{g})}{\vol(C(A,1/4),\tilde{g})}= 2^{4}.
$$
Since
$$
\frac{\vol(C(A,r_k),g')}{\vol(C(A,r_k/2),g')}\rightarrow\frac{\vol(C(A,1),\tilde{g})}{\vol(C(A,1/2),\tilde{g})}$$
and
$$\frac{\vol(C(A,r_k/2),g')}{\vol(C(A,r_k/4),g')}\rightarrow\frac{\vol(C(A,1/2),\tilde{g})}{\vol(C(A,1/4),\tilde{g})},$$
we get \eqref{3circles1} for $r=r_k$. A contradiction appears.

When $b=0$ and $a\neq 0$, we have
$$
\frac{\vol(C(A,1),\tilde{g})}{\vol(C(A,1/2),\tilde{g})}=
\frac{\vol(C(A,1/2),\tilde{g})}{\vol(C(A,1/4),\tilde{g})}=\frac{1}{2^4}.
$$
We can get another contradiction by the same argument.

To prove the lemma, now we only need to show \eqref{3circles2} does not hold. When \eqref{3circles2} holds, we can pick
$r_0$ such that
$$
\vol(C(A,2^{-k}r_0),g')>2^3\vol(C(A,2^{-k+1}r_0),g'),
$$
which contradicts $\vol(B_1,g')<+\infty$.
\endproof

Using the same method, or applying Klein transformation, we have
the following:

\begin{lem}\label{sing2}
Let $u\in W^{3,2}_{loc(\R^4\setminus B_R)}$ and $g'=u^2g_{\R^4}$.  If
$$
\vol(\R^4\setminus B_R,g')+\int_{\R^4\setminus B_R}(|K(g')|^2+|\nabla_{g'}K(g')|^2)dV_{g'}<+\infty,
$$
then, when $r$ is large enough, there holds true
$$
\vol(C(A,r))<\frac{1}{2^3} \vol(C(A,r/2)).
$$
\end{lem}

Now, we are in the position to prove the main theorem of this section:

\begin{thm}\label{kcurvature}
Let $(M,g)$ be a closed $4$-dimensional Riemannian manifold with constant scalar curvature.  Let $u_{k}\in W^{3,2}(M,g)$ be a positive function and $g_{k}=u_{k}^{2}g$. Assume
$$
vol(M,g_{k})=a_{0}\s\s\mbox{and}\s\s \int_{M}(|\nabla_{g_{k}}K(g_{k})|^{2}+|K(g_{k})|^{2})dV_{g_{k}}<\Lambda,
$$
where $a_{0}>0$ and $\Lambda >0$. Then,

1) as $(M,g)$ is not conformal to $\mathbb{S}^{4}$, $u_{k}$ converges in $W^{3,2}(M,g)$ to a positive function weakly.

2) as $M=\mathbb{S}^{4}$, there exist  M\"obius transformation $\sigma_k$ such that $\sigma_k^*(g_k)$ converges to
$W^{3,2}$-metric weakly in $W^{3,2}$.
\end{thm}

\proof

After passing to a subsequence, we find a finite  set $\CS$ such that
$$\lim_{r\rightarrow 0}\liminf_{k\rightarrow \infty}\int_{B_{r}(x)}R_{k}^{2}u_{k}^{4}>\frac{\tau}{2},\s\s x\in \CS
$$
and
$$
\lim_{r\rightarrow 0}\limsup_{k\rightarrow \infty}\int_{B_{r}(x)}R_{k}^{2}u_{k}^{4}<\frac{\tau}{2},\s\s x\notin \CS.
$$
For more details we refer to Section 5 in \cite{Li-Zhou}.

By Lemma \ref{4dim.convergence}-\ref{epsilon}, and Corollary  \ref{case1}, we can find $c_{k}>0$ such that $c_{k}u_{k}$ converges to a positive function $\phi$ weakly in $W^{3,2}_{loc}(M\setminus\CS)$.
When $\CS=\emptyset$, $c_ku_k$ converges weakly in $W^{3,2}(M,g)$, then
it follows from $\vol(M,g)=a_0$ that a subsequence of $\{c_k\}$ converges to a positive constant. Hence $\CS=\emptyset$ implies that $u_k$
converges weakly in $W^{3,2}(M,g)$.\medskip

Now, we assume $\CS\neq\emptyset$. First, we consider the case $M$ is not conformal to $\mathbb{S}^4$.
For this case, we claim  that
$$\int_{M\setminus\CS}\phi^{4}dV_{g}<\infty.$$
Assume this is not true. Since
$$
\lim\limits_{k\rightarrow \infty}\int_{M\setminus\CS}(c_{k}u_{k})^{4}dV_{g}=\int_{M\setminus{\CS}}\phi^{4}dV_{g},
$$
it follows from  $\int_{M}u_{k}^{4}=a_{0}$ that $c_{k}\rightarrow \infty$.\medskip

Let $\hat g_{k}=c_{k}^{2}g_{k}=c_{k}^{2}u_{k}^{2}g$. We have
$$
\int_{M}|K(\hat g_{k})|^{2}dV_{\hat g_{k}}=\int_{M}|K(g_{k})|^{2}dV_{g_{k}}\leq\Lambda,
$$
and
$$
\int_{M}|\nabla_{\hat g_{k}}K(\hat g_{k})|^{2}dV_{\hat g_{k}}=\frac{1}{c_{k}^{2}}\int_{M}|\nabla_{g_{k}}K(g_{k})|^{2}dV_{g_{k}}\rightarrow 0.
$$
Then, $K(\hat{g}_k)$ converges to a constant weakly in $W^{1,2}_{loc}(M\setminus\CS)$. Noting that
$$\int_{M}|K(\phi)|^2\phi^4dV_g\leq \Lambda,$$
we get $K(\phi)= 0$, which implies that $R(\phi)=0$.

By Corollary \ref{singularity}, we know that $(M,\phi^{2}g)$ is complete. On the other hand, each end of $M\setminus\CS$ is collared topologically by
$S^3\times\R$. Therefore, we conclude that $(M\setminus\mathcal{S},\phi^2g)$ is just $\R^4$. This contradicts the
assumption that $M$ is not conformal to $\mathbb{S}^4$. Thus, we get the claim.\medskip

Choose a normal chart of a point $p\in\CS$. By the definition of $\CS$, we can get a sequence $(x_{k},r_{k})$, such that $x_k\rightarrow 0$ and $r_k\rightarrow 0$ and
$$
\int_{B_{r_k}(x_{k})}|R(g_{k})|^{2}dV_{g_{k}}=\frac{\tau}{2},
$$
$$
\int_{B_{r}(y)}|R(g_{k})|^{2}dV_{g_{k}}\leq\frac{\tau}{2}, \s\s
\forall y\in B_\delta(0),\s r\leq r_k.
$$
Let $v_{k}(x)=r_{k}u_{k}(x_{k}+r_{k}x)$ and $$g_{k}'=r_{k}^{2}u_{k}^{2}(x_{k}+r_{k}x)g(x_{k}+r_{k}x).$$
It is easy to check that
$$\|K(g_{k}')\|_{W^{1,2}(B_R,g_{k}')}<C(R),\s\s \forall R,$$

$$\int_{B_1}|R(g_k')|^2dV_{g_k'}=\frac{\tau}{2},\s\mbox{and}\s
\int_{B_1(y)}|R(g_k')|^2dV_{g_k'}\leq\frac{\tau}{2},\s\s \forall y.
$$
By Lemma \ref{4dim.convergence}-\ref{epsilon}, Lemma \ref{case1} and Lemma \ref{example},
there exists a sequence of positive numbers $\{c_{k}'\}$ such that $c_{k}'v_{k}$ converges weakly to a positive function $\psi$ in $W^{3,2}_{loc}(R^{4})$ weakly. Noting $\int_{B_R}v_k^4<a_0$, we have $\inf_kc_{k}'>0$.\medskip

We claim that $$\int_{\R^4}\psi^4dx=\vol(\R^4,\psi^2g_{\R^4})<+\infty.$$
Assume this is not true. By a similar argument with the proof of $\int_M\phi^4<+\infty$,
we can get $c_k\rightarrow+\infty$ and $K(\psi)=0$.
Noting that
$$
\int_{B_1}|R(g_k')|^2dV_{g_k'}=\int_{B_1}|R({c_k'}^2g_k')|^2
dV_{c_k^2g_k'}=\frac{\tau}{2}
$$
we get $$\int_{B_1}|R(\psi)|^2dV_{\psi^2g_{\R^n}}=\frac{\tau}{2},$$
which is impossible. Therefore, the claim is  true.

Let $A'$ be an open ball in $\mathbb{S}^{n-1}$ such that, after passing to a subsequence,
$$
\int_{x_k+A'_\delta}|R(g_k)|^2dV_{g_k}<\frac{\tau}{2}.
$$
Let $A\subset  A'$ be a closed ball  in $\mathbb{S}^{n-1}$, and $\delta$
be sufficiently small.
Take $t_{k}\in [\frac{r_{k}}{\delta},\delta]$, such that
$$
\vol(C(A,t_k)+x_k,g_k)=\inf_{t\in[\frac{r_{k}}{r},r]}\vol(C(A,t)+x_k,{g_{k}}).
$$
By Lemma \ref{sing1}, for any fixed sufficiently small $r$, we have
$$
\frac{\vol(C(A,r)+x_k,g_k)}
{\vol(C(A,r/2)+x_k,g_k)}\rightarrow \frac{\vol(C(A,r),\phi^2 g)}
{\vol(C(A,r/2),\phi^2 g)}> 2^3.
$$
Then,  $t_k\rightarrow 0$. By the same argument, we deduced from Lemma \ref{sing2} that $\frac{t_k}{r_k}\rightarrow+\infty$.\medskip

Set
$$
\tilde v_{k}=t_{k}u_{k}(x_{k}+t_{k}x),\s\s \tilde g_{k}=\tilde v_{k}g(x_{k}+t_{k}x).
$$
Using the same method as we get $\phi$, we can find a finite set $\tilde\CS$ and a number $\tilde{c}_k$, such that
$\tilde{c}_k\tilde{v}_k$ converges to a positive function $v$ weakly in
$W^{3,2}_{loc}(\R^4\setminus(\{0\}\cup\CS))$. By the definition of $A$, we have
$$\CS\cap \{tA:t>0\}=\emptyset,$$
hence it follows
\begin{equation}\label{thin}
\vol(C(A,1),v^2g_{\R^n})=\inf_{t>0}\vol(C(A,t),{v^2g_{\R^n}}).
\end{equation}
Then, by the same arguments as we derive $\int_M\phi^4<+\infty$, we also obtain that
$\tilde{c}_k\rightarrow+\infty$ and $K(v)=0$. Then $v$ is a positive
harmonic function defined on $\R^4\setminus(\CS\cup\{0\})$. Furthermore,
by Theorem 9.8 in \cite{Sheldon-Paul-Wade}, for any $x_0\in \CS\cup\{0\}$ we have
 $$v(x)\sim c(x_0)|x-x_0|^{-2},$$
where $c(x_0)$ is a nonnegative constant.

Let $$\CS'=\{x\in\CS\cup\{0\}: c(x)>0\}.$$
If $\CS'$ is nonempty, then, $(\R^4\setminus\CS',v^2g_{\R ^4})$ is a complete flat manifold, whose ends are collared topologically by $S^3\times \R$. It is impossible. This means that $\CS'=\emptyset$, hence $v\in C^\infty(\R^4)$ which contradicts \eqref{thin}. Therefore, we finish the proof of 1).

Next, we consider the case $(M,g)$ is conformal to $\mathbb{S}^{4}$.
Let $P$ be the  stereographic projection from $\mathbb{S}^{4}$ to $\R^{4}$, which sends
$x_0\in\CS$ to $0\in\R^4$. Under the coordinate system defined by $P$,
as before, we can find $x_k\rightarrow 0$, $r_k\rightarrow 0$, and $c_k'$,
such that $c_k'r_ku_k(x_k+r_kx)$ converges to a positive function $\psi$, which
satisfies
$$
\int_{\R^4}\psi^4dx<+\infty.
$$
Let $\sigma_{k}(y)=P^{-1}(r_{k}P(y)+x_k)$.
It is well-known that $\sigma_k$ defines a M\"obius transformation of $\mathbb{S}^4$.
It is easy to check that for the new sequence $g_k'=\sigma_k^*(g_k)=(u_k')^2g_{\mathbb{S}^4}$, there exist
$c_k$ and a finite set $\CS'$, such that $c_ku_k$ converges weakly
in $W^{3,2}(M\setminus\CS')$
to a positive function $\phi$, which satisfies $\int\phi^4<+\infty$. Then, following the
arguments taken in 1), we complete the proof easily.
\endproof

\section{H\'elein's convergence Theorem }
The arguments in the
previous sections seem useless to the Gauss equation in 2 dimensional case under Gauss curvature condition. However, we can apply them to study the convergence of a $W^{2,2}$-conformal immersion with bounded $\|A\|_{L^2}$ to give a generalized H\'elein's Convergence Theorem.

In \cite{Kuwert-Li}, we defined the $W^{2,2}$-conformal immersion  as follows:
\begin{defi}\label{defconformalimmersion}
Let $(\Sigma,g)$ be a Riemann surface. A map $f\in W^{2,2}(\Sigma,g,\mathbb{R}^n)$
is called a conformal immersion, if the induced metric
$g_{f} = df\otimes df$ is given by
$$
g_{f} = e^{2u} g \quad \mbox{ where } u \in L^\infty(\Sigma).
$$
For a Riemann surface $\Sigma$ the set of all
$W^{2,2}$-conformal immersions  is denoted by
$W^{2,2}_{conf}(\Sigma,g,\R^n)$. When $f\in W^{2,2}_{loc}
(\Sigma,g,\R^n)$ and $u\in L^\infty_{loc}(\Sigma)$, we say
$f\in W^{2,2}_{conf,loc}(\Sigma,g,\R^n)$.
\end{defi}

H\'elein's Convergence Theorem was first proved by H\'elein \cite{Helein}. An optimal version of the theorem was stated in  \cite{Kuwert-Li} as follows:
\begin{thm}\label{Helein} Let $f_k\in W^{2,2}(D,\R^n)$
be a sequence of conformal immersions with induced metrics
${(g_{f_{k}})}_{ij} = e^{2u_k} \delta_{ij}$ and satisfy
\begin{equation}\label{ThmHelein}
\int_D |A_{f_k}|^2\,d\mu_{f_k} \leq \gamma <
\gamma_n =
\begin{cases}
8\pi & \mbox{ for } n = 3,\\
4\pi & \mbox{ for }n \geq 4.
\end{cases}
\end{equation}
If $\mu_{f_k}(D) \leq C$ and $f_k(0) = 0$, where $\mu_{f_k}$
is the measure defined by $f_k$, then $f_k$ is bounded in $W^{2,2}_{loc}(D,\R^n)$, and there is a subsequence such that one of the following two alternatives
holds:
\begin{itemize}
\item[{\rm (a)}] $u_k$ is bounded  and
$f_k$ converges weakly in $W^{2,2}_{loc}(D,\R^n)$ to a conformal
immersion $f \in W^{2,2}_{loc}(D,\R^n)$.
\item[{\rm (b)}] $u_k \to - \infty$ and $f_k \to 0$ locally uniformly on $D$.
\end{itemize}
\end{thm}

Note that in case of (a), $\|u_k\|_{W^{1,2}}<C$ follows from the boundness of
$\|u_k\|_{L^\infty}$ and $\|f\|_{W^{2,2}}$.

H\'elein's convergence Theorem is a very powerful tool to study variational
problem concerning Willmore functional \cite{Kuwert-Li,Riviere}. However, Theorem \ref{Helein} can not get rid of a collapsing sequence.
For this case, generally it is not true that $f_k$ converges to a non-trivial map after rescaling. For example, if $f_k=a_ke^{kz}$, which is a sequence of conformal maps from $D$ to $\C$, where $a_k$ is chosen  such that $\mu_{f_k}(D)=1$, then  $f_k$ converges to point, and for any $c_k$, $c_kf_k$ does not converge. However, in \cite{Li1} (also see \cite{Li2}) Y. Li showed that, if $f_k(D)$ can be extended to a closed surface immersed in $\R^n$ with $\|A\|_{L^2}<C$, then we can find $c_k$, such that $c_kf_k$ converges weakly in $W^{2,2}(D_r)$ for any $r$ to a conformal immersion. The proof provided in \cite{Li1} is based on the conformal invariant of Willmore functional and Simon's monotonicity formula.

In this section, we will use the John-Nirenberg inequality to give a new sufficient condition to guarantee the above assertion is still valid.

We define
$$
\rho(u_k,x)=\sup\left\{t:
\int_{D_r(x)\cap D}|\nabla u_k|^2<\epsilon_0^2\right\}.
$$
We first prove the following:

\begin{lem} Let $f\in W_{conf}^{2,2}(D,\R^{n})$. Suppose that there exists a positive number $\beta$ such that, for any $y\in \R^n$ and $r>0$,
\begin{equation}\label{measure.d}
\frac{\mu_{f}(f^{-1}(B_r(y)))}{\pi r^2}<\beta.
\end{equation}
Then there exists $\epsilon>0$ and $a>0$ such that, if $\int_D|A|^2<\epsilon$, then
$$
\inf_{x\in D_\frac{1}{4}}\rho(u,x)>a.
$$
\end{lem}

\proof If this is not true, then, we can find  a sequence of $f_k$, such that
$\int_D|A_k|^2\rightarrow 0$ and $\inf_{D_\frac{1}{4}}\rho(u_k,x)
\rightarrow 0$.
Take $x_k\in D_\frac{1}{4}$, such that $\rho(u_k,x_k)\rightarrow 0$
and $x_k\rightarrow x_0$.

Put $z_k\in D_\frac{1}{2}$ such that
$$
\frac{\rho(u_k,z_k)}{1/2-|z_k|}=\inf_{x\in D_\frac{1}{2}}\frac{\rho(u_k,x)}{1/2-|x|}:=\lambda_k.
$$
As in the proof of Corollary \ref{radius}, we have  $\rho_k:=\rho(u_k,z_k)
\rightarrow 0$, $D_{R\rho_k}(z_k)\subset D_\frac{1}{2}$ and
$$
\frac{\rho(u_k,z)}{\rho(u_k,z_k)}>\frac{1}{2},\s\s \forall z\in D_{R\rho_k}(z_k),
$$
when $k$ is sufficiently large.

Assume $z_k\rightarrow z_0$ and put
 $f_k'(x)=c_k(f_k(z_k+\rho_kz)-f(z_k))$, where $c_k$ is chosen such that
$$
\int_Du_k'=0.
$$
It is easy to see that $f_{k}'$ also satisfies \eqref{measure.d}. Then, as in proof of Corollary \ref{radius}, we have $\int_{D_R}e^{2u_k'}<C(R)$ for any $R$. Since $\int_Du_k'$ does not converge to $-\infty$, by Theorem \ref{Helein}, we know that $f_k'$ converges weakly in $W^{2,2}_{loc}(\C,\R^n)$
to an $f'\in W^{2,2}_{loc}(\C,\R^n)$ with $A_{f'}=0$. Since $f'$ is conformal, it is a  holomorphic
immersion from $\C$ to a plain $L$ in $\R^n$.

Moreover, from $$-\Delta u_k=K_{f_k}e^{2u_k}$$ we deduce that $u'$ is a harmonic function on $\R^2$ and hence $\nabla u'$
is harmonic, since $K_{f_k}e^{2u_k}$ converges to 0 in $L^1$ and
$u_k'$ converges to $u'$ weakly in $W^{1,2}_{loc}(\C)$. Obviously, we also have that for any $x\in \R^2$
$$\int_{D_\frac{1}{2}(x)}|\nabla u'|^2dx\leq\epsilon_0^2,$$
which follows from that for any $x$
$$\int_{D_\frac{1}{2}(x)}|\nabla u_k'|^2dx\leq\epsilon_0^2.$$
By mean value theorem, $\nabla u'$ is bounded. Therefore, $\nabla u'$ is a constant vector. Choosing an appropriate coordinates of $L$, we may write $f'$ as  $f'=az$ or $e^{az+b}$, where $a\neq 0$.\medskip

When $f'=az$, $u'$ is a constant.
Note that i) of Theorem \ref{Helein} implies that for any $r$
$$\|u_k'\|_{W^{1,2}(D_r)}<C(r).$$
Without loss of generality, we assume $u_k'$ converges to $u'$ weakly in $W^{1,2}_{loc}(\C)$.
Given an positive cut-off function $\eta$ which is 1 on $D_1$,
we have
\begin{eqnarray*}
\epsilon_0^2&\leq&\int_{\C}\eta|\nabla u_k'|^2=
\int_{\C}\nabla(\eta u_k')\nabla u_k'-u_k'\nabla\eta\nabla u_k'\\
&=&
\int_{\C}\eta u_k' Ke^{2u_k'}-\int_{\C}(u_k'-u')\nabla\eta\nabla u_k'
-u'\int_{\C}\nabla\eta\nabla u_k'\\
&\rightarrow& 0.
\end{eqnarray*}
This is a contradiction.\medskip

When $f'(z)=e^{az+b}$,  there exists $P_0\in L$, such that
${f'}^{-1}(\{P_0\})$ contains infinity many points. Let $m>\beta+1$.
Take
$z_1$, $\cdots$, $z_m\in {f'}^{-1}(\{P_0\})$ and choose $r>0$ and $r'>0$ such that $B_{r'}(P_0)\cap L
\subset f(D_r(z_i))$ and $f$ is injective on $D_r(z_i)$. Then we get
$$
\frac{\mu_{f'}\left({f'}^{-1}(B_{r'}(P_0))\right)}{\pi {r'}^2}=m,
$$
and hence
$$
\frac{\mu_{f_{k}'}\left({f_{k}'}^{-1}(B_{r'}(P_0))\right)}{\pi {r'}^2}>m-1>\beta,
$$
when $k$ is sufficiently large. This contradicts \eqref{measure.d}.
\endproof

\begin{thm}\label{Helein.g}
Let $f_k\in W^{2,2}_{conf}(D,\R^n)$ and satisfy \eqref{measure.d}.
Then, there exists an $\epsilon>0$ such that, if
$\int_D|A_{f_k}|^2<\epsilon$, then there exist $c_k$ such that
$c_kf_k$ converges weakly in $W^{2,2}_{loc}(D_r)$ to an
$f\in W^{2,2}_{conf,loc}(D,\R^n)$ for any $r<1$.\\
\end{thm}

\proof We only need to prove that, there
exists $c_k$, such that
$$
\int_{D_r}e^{2|u_k+\log c_k|}<C(r)
$$
for any
$r$. {
The proof goes almost the same as in the proof of Lemma \ref{example}, we omit it.
\endproof

When $f_k$ can be extended to a closed immersed surface with
$\|A\|_{L^2}<C$, by (1.3) in \cite{Simon}, we know that \eqref{measure.d} must hold true.

\begin{cor}
Let $f_k\in W^{2,2}_{conf}(D,\R^n)$, which satisfies \eqref{measure.d}. If $f_k$ satisfies
$$\int_D|A_{f_k}|^2<\gamma_n-\tau,$$
then there exist $c_k$ such that $\{c_kf_k\}$ converges weakly in $W^{2,2}_{loc}(D_r)$ to an $f\in W^{2,2}_{conf,loc}(D,\R^n)$ for any $r<1$.
\end{cor}

\proof Let $\epsilon$ be the same as in the Theorem \ref{Helein.g}.
Take $m$ such that $\frac{8\pi}{m}\cdot 5<\epsilon$. For convenience, we set $r\in(\frac{3}{4},\,1)$ and
$l=\frac{1-r}{m}$.
After passing to a subsequence, there exists $2\leq i\leq m-2$ such that
$$
\int_{D_{r+(i+2)l}\setminus D_{r+(i-2)l}}|A_{f_k}|^2<\epsilon,\s\s \forall k.
$$
By Theorem \ref{Helein}, Theorem \ref{Helein.g} and a covering argument, we know there exists $c_{k}'$ such that $c_{k}'f_{k}$ converges weakly in $W^{2,2}_{loc}$ to a function
$$f_{0}\in W^{2,2}_{conf,loc}(D_{r+(i+1)l}\setminus D_{r+(i-1)l},\R^n),$$
and
$$
\|u_{k}+\log c_{k}'\|_{L^{\infty}(D_{r+(i+1)l}\setminus D_{r+(i-1)l})}<C.
$$
In particular, there holds true
$$
\| u_{k}+\log c_{k}'\|_{L^{\infty}(\partial D_{r+il})}<C.
$$

Since $$\int_D|A_{f_k}|^2<\gamma_n-\tau,$$
by Corollary 2.4 in \cite{Kuwert-Li}, we know there exists a function $\upsilon_{k}: \C\rightarrow \R$ solving the equation $$-\Delta \upsilon_{k}=K_{c_{k}'f_{k}}e^{2(u_{k}+\log c_{k}')}$$
in $D$ and satisfying the following estimates:
$$\|\upsilon_{k}\|_{L^{\infty}(D)}\leq C.$$
The maximal principle yields $$\| u_{k}+\log c_{k}'-\upsilon_{k}\|_{L^{\infty}{(D_{r+il})}}<C,$$
hence it follows
$$
\| u_{k}+\log c_{k}'\|_{L^{\infty}{(D_{r+il})}}<C.
$$
By the fact $f_k$ satisfies the equation $$\Delta c_{k}'f_{k}=e^{2(u_{k}+\log{c_{k}})}H_{c_{k}'f_{k}}$$ for every $k$, we obtain
$$
\int_{D_{r+il}}|\Delta c_{k}'f_{k}|^{2}dx\leq e^{2\| u_{k}+\log c_{k}'\|_{L^{\infty}(D_{r+il})}}\int_{D_{r+il}}|H_{c_{k}'f_{k}}|^{2}d\mu_{c_{k}'f_{k}}\leq C.
$$
This implies $$\|c_{k}'f_{k}\|_{W^{2,2}{(D_{r+il})}}<C.$$
Thus, there exists a subsequence of $\{c_{k}'f_{k}\}$ converges weakly to a $W^{2,2}$ conformal immersion in $D_{r}$.

Applying Theorem \ref{Helein} again, we get
$$
\|\nabla (u_k+\log c_k')\|_{L^2(D_r)}+\|u_k+\log c_k'\|_{L^\infty(D_r)}<C(r).
$$
Let $$\log c_k= -\frac{1}{|D_\frac{1}{2}|}\int_{D_\frac{1}{2}}u_k.$$
By Poincar\'e inequality, we have
$$\|u_k  + \log c_k\|_{L^2(D_r)}<C.$$
Hence, it follows that $$|\log c_k-\log c_k'|<C.$$
Thus, after passing to a subsequence, $c_kf_k$ converges weakly in $W^{2,2}(D_r)$ to a conformal map.
\endproof


\begin{thebibliography}{2}

\bibitem{Sheldon-Paul-Wade}S. Axler, P. Bourdon, W. Ramey:
Harmonic function theory. Second edition. Graduate Texts in Mathematics, 137. Springer-Verlag, New York, 2001.


\bibitem{Chang-Yang1}S.Y.A. Chang, P. Yang: {Compactness of isospectral conformal metrics on ${\CS}^{3}$},
{Comment. Math. Helvetici}, {\bf 64} (1989), {no. 3}, 363-374.

\bibitem{Chang-Yang2}S.Y.A. Chang, P. Yang: Isospectral Conformal Metrics on 3-Manifolds. {\em J. Amer. Math. Soc.} {\bf 3} (1990), {no. 1}, 117-145.

\bibitem{Gursky}M. Gursky: Compactness of conformal metrics with integral bounds on curvature. {\em Duke Math. J.} , {\bf 72} (1993), {no. 2}, 339-367.

\bibitem{Helein} F. H\'elein: Harmonic maps,
conservation laws and moving frames.
Translated from the 1996 French original.
With a foreword by James Eells. Second edition.
Cambridge Tracts in Mathematics, 150.
Cambridge University Press, Cambridge, 2002

\bibitem{Kuwert-Li} E. Kuwert and Y. Li: $W^{2,2}$-conformal immersions of a closed Riemann
surface into $\R^n$.  {\em Comm. Anal. Geom.} {\bf 20} (2012),  313-340.

\bibitem{Leckband}
M. Leckband:
Moser's inequality on the ball $B^n$ for functions with mean value zero.
{\em Comm. Pure Appl. Math.}  {\bf 58}  (2005),  no. 6, 789-798.

\bibitem{Li1} Y. Li:
Weak limit of an immersed surface sequence with bounded Willmore functional. arXiv:1109.1472.


\bibitem{Li2} Y. Li: Some remarks on Willmore surfaces embedded in $\R^3$. {\em J. Geom. Anal.}  {\bf 26}  (2016),  no. 3, 2411-2424.


\bibitem{Li-Zhou} Y. Li, Z. Zhou: Conformal metric sequences
with  integral-bounded scalar curvature.  arXiv:1706.03919.

\bibitem{Riviere} T. Rivi\'ere: Variational principles for immersed surfaces with $L^2$-bounded second fundamental form. {\em J. Reine Angew. Math.}  {\bf 695}  (2014), 41-98.


\bibitem{Simon} L. Simon: Existence of surfaces minimizing
the Willmore functional, {\em Comm. Anal. Geom.}
{\bf 1} (1993),
281--326.



\end{thebibliography}
\end{document}